\documentclass{article}
\usepackage[top=2cm, bottom=2cm, left=3cm, right=3cm]{geometry}
\usepackage{authblk}
\usepackage[utf8]{inputenc}
\usepackage[T1]{fontenc}
\usepackage{cite}
\usepackage{amsmath,amsthm, amssymb}
\usepackage[english]{babel}
\usepackage{enumitem}
\usepackage{mathtools}
\usepackage{subfiles}
\usepackage{hyperref}
\usepackage[backgroundcolor=white, textsize=tiny]{todonotes}

\theoremstyle{plain}
\newtheorem{thm}{Theorem}[section]
\newtheorem{lem}[thm]{Lemma}
\newtheorem{obs}[thm]{Observation}
\newtheorem{prop}[thm]{Proposition}
\newtheorem{cor}[thm]{Corollary}

\theoremstyle{definition}
\newtheorem{defn}[thm]{Definition}
\newtheorem*{defn*}{Definition}
\newtheorem{exmp}[thm]{Example}

\theoremstyle{remark}
\newtheorem*{rem}{Remark}

\newcommand{\diff}{{d}}
\newcommand{\omegabar}{\bar{\omega}}
\newcommand{\ld}{{L}}

\begin{document}
\baselineskip 12pt

\title{Density-valued symplectic forms from a multisymplectic viewpoint}

\author[1,2]{Laura Leski}
\author[2,3]{Leonid Ryvkin}
\affil[1]{University of Würzburg}
\affil[2]{University of Göttingen}
\affil[3]{Universit\'e Claude-Bernard Lyon 1}

\maketitle
\begin{abstract}
	We give an intrinsic characterization of multisymplectic manifolds that have the linear type of density-valued symplectic forms in each tangent space, prove Darboux-type theorems for these forms, and investigate their symmetries.
\end{abstract}
\tableofcontents

\section{Introduction}
Symplectic geometry lies at the core of the Hamiltonian formulation of classical mechanics. Darboux theorem, asserting that locally any symplectic form is equivalent to the standard symplectic form on the cotangent bundle $T^*\mathbb R^n$ of Euclidean space, is a cornerstone of symplectic geometry: It greatly simplifies local computations, and implies that symplectic manifolds have no local invariants (other than the dimension) - opening the door for the question about global invariants - i.e., symplectic topology. In this article, we are going to study the analogue of the Darboux theorem for a specific class of multisymplectic forms, i.e., non-degenerate closed differential forms of degree higher than two.\\

Multisymplectic geometry was initially introduced to extend the symplectic perspective from particle mechanics to classical field theories \cite{MR0334772}, cf. \cite{MR2882772,ww} for more recent accounts including some historical remarks. Multisymplectic forms have a certain number of new phenomena, not present in the symplectic case (cf. \cite{RYVKIN20199} for an introduction to them):

\begin{enumerate}
	\item Fixing a dimension $n$ and a form rank $k$, there are, typically, uncountably many inequivalent linear models of non-degenerate alternating forms \cite[Proposition 6.3]{martinet70}.
	      On connected manifolds, multisymplectic forms with a non-constant linear model can be constructed.
	\item Two multisymplectic forms with identical, constant linear models need not be equivalent at the level of germs. First results of this flavor go back to \cite{MR801210}, with the case of classical field theory being first treated in \cite{martin} (cf. also \cite{sevwurz} for a more detailed exposition including a more extensive historical account). This means that having a Darboux-type theorem can require additional conditions beyond the closedness of the multisymplectic form.
	\item The notion of Hamiltonian vector field still exists. It describes a vector field whose contraction into the multisymplectic form is exact. However, in the general multisymplectic setting, this notion is weaker: There might be very few such fields, and their corresponding differential forms (which in the symplectic case would be just ordinary functions) do not carry a Lie (or Poisson) algebra structure. However, they can be turned into a Leibniz or $L_\infty$-algebra structure \cite{MR2892558}, which allows the use of methods from homotopical algebra in the multisymplectic setting.
\end{enumerate}

Hence, a very natural procedure for studying multisymplectic manifolds from a geometric viewpoint is fixing a linear normal form (item 1 above), describing its obstructions to a Darboux theorem (item 2), and studying its symmetries (item 3). In this note, we follow this procedure for differential forms that (on the level of an individual tangent space) look like the external wedge product of a symplectic form with a volume form, i.e.
\begin{align}\label{eq:intro}
	\omega= \left(\sum_{i=1}^m \diff x_{2i-1}\wedge \diff x_{2i}\right)\wedge \diff y_1 \wedge \hdots \wedge \diff y_k .
\end{align}

We call such forms dvs-forms (density-valued symplectic forms), to underline one of their motivations: Let $(P,\eta)$ be a symplectic manifold and $(\Sigma,\sigma)$ a compact manifold equipped with a volume form. We can obtain a (pre-)symplectic structure on the mapping space $\mathcal E=C^\infty(\Sigma, P)$ (seen as an infinite-dimensional manifold \cite[Theorem 42.1]{krieglMichor}), by transgressing the form $\eta\wedge\sigma$ to it. More generally, a (non-trivial) $P$-bundle $E$ over $\Sigma$ equipped with a form which locally looks like \eqref{eq:intro} would yield a pre-symplectic structure on $\Gamma(E)$, the section space of $E\to \Sigma$.\\

The article is structured as follows: In Section \ref{sec:darboux}, we establish a Darboux-type theorem for dvs-forms, which necessitates an additional integrability condition beyond the closedness of the multisymplectic form at hand. We then proceed to show that when there are sufficiently many symplectic directions, this condition is automatically satisfied. In Section \ref{sec:symmetry}, we then study the symmetries of dvs-forms and relate them to conformally Hamiltonian symmetries of symplectic manifolds. Finally, in Section \ref{sec:cosymplectic} we discuss how cosymplectic manifolds and other structures stemming from Poisson geometry relate to dvs-forms.\\

The main results of the article are Theorem \ref{thm:ourdarboux} and Theorem \ref{thm:autoflat}, which generalize results of \cite{MR801210} on infinitesimally transitive three-forms.

\subsubsection*{Acknowledgements}
The authors would like to thank Christian Blohmann, Alfonso Garmendia, Aliaksandr Hancharuk, Frédéric Hélein, Michele Schiavina, Thomas Schick, Aldo Witte, and Tilmann Wurzbacher for valuable discussions related to this work. L.R. was supported by the  DFG grant Higher Lie Theory - Project number 539126009. For the purpose of Open Access, a CC-BY-NC-SA public copyright licence has been applied by the authors to the present document and will be applied to all subsequent versions up to the Author Accepted Manuscript arising from this submission. Parts of this article are based on the Bachelor's thesis \cite{Leski2022} for the first author. We would like to thank the referees for their careful reading of the manuscript and suggestions on how to improve its exposition.

\section{A Darboux theorem}
\label{sec:darboux}

\subsection{Flatness of multisymplectic forms}
\begin{defn}A multisymplectic (or $k$-plectic) manifold is a couple $(M,\omega)$, where $M$ is a smooth manifold and $\omega\in \Omega^{k+1}(M)$ differential form on $M$, satisfying:
	\begin{enumerate}
		\item closedness: $d\omega=0$,
		\item non-degeneracy: the map $v\mapsto \iota_v\omega:TM\to \Lambda^{k}T^*M$ is injective.
	\end{enumerate}
\end{defn}

Non-degeneracy is a fiberwise condition and can be defined for alternating forms $\eta\in\Lambda^{k+1}V^*$ for any vector space $V$.

\begin{exmp}Symplectic manifolds are precisely 1-plectic manifolds and volume forms on a manifold $M$ correspond to $(\mathrm{dim}(M)-1)$-plectic structures on $M$.
\end{exmp}

The classical Moser theorem asserts that any volume form on a $k$-dimensional manifold can be locally written as $dy_1\wedge \ldots \wedge dy_k$ upon choosing appropriate coordinates $(y_1,\ldots ,y_k)$. Similarly, the Darboux theorem asserts that locally, any symplectic form (in dimension $2m$) looks like $dx_1\wedge dx_2 + \ldots  +dx_{2m-1}\wedge dx_{2m}$. Both theorems can be proven in two steps. The first step is linear algebra: One shows that all non-degenerate alternating $n$-forms in $\mathbb R^n$ (resp. 2-forms in $\mathbb R^{2m}$) can be transformed into one another by a base change. The second step is so-called flatness:

\begin{defn}A differential form $\omega$ is called flat near $p$ if there are local coordinates near $p$ in which $\omega$ has constant coefficients. A form is called flat if it is flat near all points.
\end{defn}

Unfortunately, neither part of the Moser/ Darboux theorems holds for general multisymplectic manifolds.
On the linear level, the number of 'linear normal forms' for $k$-plectic structures in dimension $n$ corresponds to the number of non-degenerate orbits of the $GL(n,\mathbb R^{n})$-action on $\Lambda^{k+1}\mathbb R^{n,*}$. For $n>8$ and $k\in\{2,\ldots ,n-4\}$ this number is infinite (cf. \cite[Proposition 6.3]{martinet70}). The only cases where there is a unique non-degenerate orbit are symplectic forms, volume forms, and 3-forms in dimension 5. However, already in this case, one can see how flatness can fail:

\begin{exmp}[\cite{MR801210}]\label{ex:nonflat3form}
	The form $(dx_1\wedge dx_2+dx_3\wedge dx_4)\wedge (dy+x_1dx_3)\in \Omega^3(\mathbb R^5)$ is multisymplectic but nowhere flat.
\end{exmp}

To see the non-flatness of the above form, \cite{MR801210} constructed the following invariant:

\begin{defn} Let\label{F(omega)undD(omega)} $\omega\in \Omega(M)$ be a differential form. We set
	$F(\omega)= \bigsqcup_{p\in M} F(\omega_p)$,
	where $F(\omega_p)=\{\alpha\in T_p^*M~|~\alpha\wedge \omega_p=0\}$.
	We also write $D(\omega_p):=F(\omega_p)^{ann}\subset T_pM$ (resp. $D(\omega):=F(\omega)^{ann}\subset TM$) for their annihilators.
\end{defn}

Since for constant coefficient forms the space $F(\omega_p)$ does not depend on $p$ one has:

\begin{obs}
	If $\omega\in \Omega(M)$ is flat, then $D(\omega)$ is an involutive distribution. (The latter is equivalent to the ideal generated by $F(\omega)$ to be differential, i.e. to the condition $d\Gamma(F(\omega))\subset \Gamma(F(\omega))\wedge \Omega^{1}(M)$).
\end{obs}

In the case of multisymplectic 3-forms in dimension 5, the converse is also true:
\begin{lem}[Theorem 1.1 in \cite{MR801210}]
	Let $M$ be 5-dimensional and $\omega\in \Omega^{3}(M)$ multisymplectic. If $D(\omega)$ is involutive, then $\omega$ is flat.
\end{lem}

In fact, \cite{MR801210} carried out a local classification (up to diffeomorphism) of all infinitesimally transitive 2-plectic structures in 5-dimensional space, including not necessarily flat ones.  \\

\subsection{Normal forms of density-valued symplectic forms}
The goal of this section is to give an intrinsic characterization of multisymplectic forms, which can be interpreted locally as the external wedge product of a symplectic form with a volume form.  We start with the pointwise, i.e., linear, description:

\begin{lem}\label{lemlin}
	Let $V$ be a $(2m+k)$-dimensional vector space and $\omega\in\Lambda^{k+2}V^*$ a non-degenerate alternating $k+2$-form. If the space $F(\omega):=\{\alpha \in V^*|\alpha\wedge\omega=0\}$  is $k$-dimensional, then $V$ admits a basis $e_1,\ldots ,e_{2m}, f_1,\ldots ,f_k$  such that
	\begin{align}\label{eq:linform}
		\omega=  (e_1^*\wedge e_2^*+ \ldots +e_{2m-1}^*\wedge e_{2m}^*)\wedge f_1^*\wedge \ldots \wedge f_k^*.
	\end{align}
\end{lem}

\begin{proof}
	The statement can be shown by extending a basis $f_1^*, \ldots , f_k^*$ of $F(\omega)$ to a basis of $V^*$ and using the fact that for any $\alpha\in V^*\backslash \{0\}$ and $\eta\in \Lambda V^*$
	$$\alpha\wedge \eta=0\Longleftrightarrow\exists\bar\eta: \eta=\bar\eta\wedge \alpha.$$
	By applying this statement multiple times, one can see that $\omega=\omegabar\wedge f_1^*\wedge \ldots \wedge f_k^*$, where $\omegabar$ is a non-degenerate 2-form on $D(\omega)$. One then applies the standard symplectic basis theorem.
\end{proof}

\begin{rem} For $k=1$ the above Lemma reproduces  \cite[Proposition 7]{vanzuraCharacterizationOneType2004}.
\end{rem}

We can now prove our central theorem:

\begin{thm}[\cite{Leski2022}]\label{thm:ourdarboux}
	Let $M$ be a $(2m+k)$-dimensional manifold and $\omega$ a closed non-degenerate $k+2$-form.  If the bundle $F(\omega)$ has (fiberwise) constant dimension $k$ and its annihilator distribution  $D(\omega)=F(\omega)^{ann}\subset TM$ is involutive, then there are local coordinates $(x_1,\ldots ,x_{2m}, y_1,\ldots ,y_k$), such that
	\begin{align}
		\label{vardarbouxform}
		\omega= \left(\sum_{i=1}^m \diff x_{2i-1}\wedge \diff x_{2i}\right)\wedge \diff y_1 \wedge \hdots \wedge \diff y_k .
	\end{align}
\end{thm}

\begin{proof} Let $p\in M$. We start by considering \(F(\omega)\) and \(D(\omega)\) defined as in Definition \ref{F(omega)undD(omega)}. As \(D(\omega)\) is involutive of dimension \(2m\), we can find local coordinates $(x_1,\ldots ,x_{2m},y_1,\ldots ,y_k)$ centered at $p$ such that $\partial_{x_i},\hdots,\partial_{x_{2m}}$ span $D(\omega)$, i.e. $dy_1,\ldots ,dy_k$ span $F(\omega)$. Thus, \(\omega\) can be rewritten as follows
	\begin{align}
		\label{yinomegaseparated}
		\omega=\omegabar\wedge \diff y_1\wedge \hdots \wedge \diff y_k ,
	\end{align}
	where \(\omegabar=\iota_{\partial_{y_1}}\hdots \iota_{\partial_{y_k}}\omega\).

	The splitting $TM=D(\omega)\oplus \mathrm{span}(\partial_{y_1},\hdots,\partial_{y_k})$ into a pair of transversal foliations induces a bigrading on the de Rham complex. The de Rham differential $d$ splits into $d=\diff^x+\diff^y$, where $\diff^x$ has bidegree $(1,0)$ and $\diff^y$ bidegree $(0,1)$. We note that $\diff^x$ and $\diff^y$ satisfy a Poincar\'e Lemma, i.e. their cohomologies are locally trivial \cite[Chapter 5]{vaisman2016cohomology}.\\

	We show that \(\omegabar\) is closed and non-degenerate with respect to the \(x_i\)-directions and use this to derive the flatness using a Moser-type argument.\\

	For the closedness, we calculate:
	\begin{align*}
		0=\diff \omega & = \diff(\omegabar\wedge \diff y_1\wedge\hdots\wedge\diff y_k)=(\diff \omegabar)\wedge\diff y_1\wedge \hdots\wedge\diff y_k                        \\
		               & =(\diff^y \omegabar +\diff^x \omegabar)\wedge \diff y_1\wedge\hdots\wedge\diff y_k= \diff^x \omegabar\wedge \diff y_1\wedge\hdots\wedge\diff y_k.
	\end{align*}
	In the last step, we used that the extra \(\diff y_i\), which we obtain by computing the differential in the \(y_i\)-directions, will vanish when wedged with the corresponding \(\diff y_i\) in the wedge product. Additionally, the resulting term can only be equal to zero if \(\diff^x \omegabar\) is equal to zero. Hence, \(\omegabar\) is a $\diff^x$-closed form.\\

	Next, we want to show that the restriction \(\omegabar|_D\) is non-degenerate. To do so, we can use the fact that \(\omega\) is non-degenerate.  So let \(X\in D\) and suppose \(\iota_X\omega\neq 0\), then we can do the following computation
	\begin{align*}
		0 & \neq \iota_X\omega=\iota_X(\omegabar\wedge (\diff y_1\wedge \hdots\wedge\diff y_k))= (\iota_X\omegabar)\wedge \diff y_1\wedge\hdots\wedge\diff y_k.
	\end{align*}
	Hence, for all \(X\in D(\omega)\backslash \{0\}\) it holds that \((\iota_X\omegabar)|_{D(\omega)}\)  can not be zero, i.e., \(\omegabar|_{D(\omega)}\) is non-degenerate.\\

	Now, we can go on to actually prove the flatness of \(\omega\). For this purpose, we evaluate \[\omegabar=\sum_{i=1}^{2m}\sum_{j>i}^{2m}f_{ij}(x,y)\diff x_i\wedge \diff x_j\] at \(0\) and get a constant coefficient form
	\begin{align*}
		\omegabar_0=\sum_{i=1}^{2m}\sum_{j>i}^{2m}f_{ij}(0)\diff x_i\wedge \diff x_j.
	\end{align*}
	Additionally,  we introduce a homotopy between \(\omegabar_0\) and \(\omegabar=:\omega_1\)
	\begin{align*}
		{\omegabar_t}=t\cdot\omegabar+(1-t)\cdot\omegabar_0
	\end{align*}
	for \(t\in [0,1]\). Accordingly, let \(\omega_t\) be the \((n+2)\)-form \(\omega_t:=\diff y_1\wedge \hdots \wedge \diff y_k \wedge \omegabar_t\). We note that \(\omega_0= \omegabar_0\wedge \diff y_1\wedge \hdots \wedge \diff y_k\) is a constant coefficient form.\\

	As \(\omegabar_t\) is \(\diff^x\)-closed its derivative with regard to \(t\), which we denote as \(\dot{\omegabar}_t=\omegabar-\omegabar_0\), is also \(\diff^x\)-closed. Now, we can use the foliated Poincaré Lemma. It says, that in an open neighborhood around \(0\), there exists a form \(\theta\) s.t. \(\diff^x \theta_t=\omegabar-\omegabar_0=\dot{\omegabar}_t\). Since $\dot{\omegabar}_t$ is a $(2,0)$-form,  \(\theta\) can be written in the form \(\theta=\sum_{i=1}^{2m} g_i(x,y)\diff x_i\), i.e., as an element of $\Omega^{1,0}$.

	Since $\omegabar_t$ is non-degenerate with respect to the $D(\omega)$-directions at the origin for all $t$, there is an open neighborhood of the origin, where \(\omegabar_t\) is non-degenerate with respect to the $D(\omega)$-directions.

	In this neighborhood, the equation \(\iota_{X_t}\omegabar_t=-\theta\) has a unique solution $X_t$, among the time-dependent vector fields with values in $D(\omega)$. \\

	Now, we show that the time-1 flow of $X_t$ intertwines \(\omega_0\) and \(\omega_1\). For that, we first show that \(\ld_{X_t}\omega_t=0\), where $\ld$ denotes the (time-dependent) Lie derivative:
	\begin{align*}
		\ld_{X_t}\omega_t & =\ld_{X_t}(\omegabar_t\wedge \diff y_1\wedge \hdots\wedge \diff y_k)                                                                                                  \\
		                  & =  \ld_{X_t}\omegabar_t \wedge\diff y_1\wedge \hdots \wedge\diff y_k                                                                                                  \\
		                  & = (\diff \iota_{X_t} \omegabar_t+\iota_{X_t}\diff \omegabar_t+\Dot{\omegabar}_t) \wedge\diff y_1\wedge \hdots \wedge\diff y_k                                         \\
		                  & =(\diff(-\theta)+\iota_{X_t}\diff \omegabar_t+\Dot{\omegabar}_t)\wedge\diff y_1\wedge \hdots \wedge\diff y_k                                                          \\
		                  & =(-\diff^x \theta-\diff^y \theta+\iota_{X_t}(\underbrace{\diff^x \omegabar_t}_{=0}+\diff^y\omegabar_t)+\Dot{\omegabar}_t)\wedge\diff y_1\wedge \hdots \wedge\diff y_k \\
		                  & =(-\dot{\omegabar}_t-\diff^y \theta+\iota_{X_t}\diff^y\omegabar_t+\Dot{\omegabar}_t)\wedge\diff y_1\wedge \hdots \wedge\diff y_k                                      \\
		                  & = (-\diff^y \theta+\iota_{X_t}\diff^y\omegabar_t)\wedge\diff y_1\wedge \hdots \wedge\diff y_k
		=0.
	\end{align*}
	In the last step, we used that the differential with respect to the \(y\)-directions adds an extra \(\diff y_i\) to the wedge product, which vanishes because there already is another \(\diff y_i\) in the wedge product. Now we further shrink the open neighborhood of the origin, such that the flow $\phi_t$ of $X_t$ is well-defined for $t\in [0,1]$ and then use the fact that  \(\ld_{X_t}\omega_t=\frac{\diff}{\diff t}\phi_t^*\omega_t\). Thus $\phi_t^*\omega_t$ is constant, in particular $\phi_1^*\omega_1=\phi_0^*\omega_0=\omega_0$.
	Hence, \(\phi_1^{-1}\) is a diffeomorphism, that transforms \(\omega=\omega_1\) to a form with constant coefficients. Using Lemma \ref{lemlin}, this constant coefficient form can even further be transformed into the form as in Equation \eqref{vardarbouxform} by a linear transformation.
\end{proof}

\begin{defn}
	We call \emph{dvs-form} a differential form, which satisfies the condition of Theorem \ref{thm:ourdarboux}, i.e., which can locally be written in the form of \eqref{vardarbouxform}. The acronym dvs comes from density-valued symplectic.
\end{defn}

\begin{exmp}[Symplectic forms] As a first edge case of the construction, we can observe that if $D(\omega)=\{0\}$, then there are no $y$ variables and the above proof exactly reconstructs Weinstein's proof of Darboux theorem using Moser's method, i.e., that any symplectic form locally looks like $dx_1\wedge dx_2 + \ldots +dx_{2m-1}\wedge dx_{2m}$
\end{exmp}

\begin{exmp}[Volume forms] Formally speaking, the case of volume forms ($m=0$ or $m=1$) can not appear in the above theorem, since all coordinates would be $y$ coordinates, i.e., there are no $x$ directions left. However, we could pick any set of local coordinates $y_1,\ldots ,y_{k}$ and relabel $y_{k-1},y_{k}$ to $x_1,x_2$ and the rest of the proof would go through, implying that any volume form locally looks like $dy_1\wedge \ldots \wedge dy_k$. Hence, we could see the theorem as a simultaneous generalization of the normal form theorems for volumes and symplectic forms.
\end{exmp}

\begin{exmp}[\(D(\omega)\) not integrable] This is an immediate generalization of Example \ref{ex:nonflat3form}. Consider the form $$\omega= (dx_1\wedge dx_2+dx_3\wedge dx_4)\wedge(dy_1+ x_2dx_4)\wedge dy_2\wedge \ldots \wedge dy_k\in \Omega^{k+2}(\mathbb R^{k+4}).$$
	It clearly is multisymplectic and pointwise of the form \eqref{eq:linform}. However, since $d(dy_1+ x_2dx_4)=dx_2\wedge dx_4$ is never a multiple of $(dy_1+ x_2dx_4),dy_2, \ldots , dy_k$, it is not flat at any point.
\end{exmp}

\subsection{Automatic involutivity}
In fact, $D(\omega)$ may be non-involutive only for $m=2$. This has been observed for infinitesimally transitive 3-forms in \cite{MR801210}. We give here a proof working in the general setting and not relying on infinitesimal transitivity:

\begin{thm}\label{thm:autoflat}
	Let $M$ be a $(2m+k)$-dimensional manifold and $\omega$ a closed non-degenerate $k+2$-form.  If the bundle $F(\omega)$ has (fiberwise) constant dimension $k$ and $m\neq 2$, then $D(\omega)$ is involutive.
\end{thm}

\begin{proof}
	Without loss of generality, we may assume $m>2$. Let $E$ be an involutive distribution such that $D(\omega)\oplus E=TM$. Let $\alpha_1,\ldots ,\alpha_k$ be a local frame for $F(\omega)$ and pick local coordinates $(x_1,\ldots ,x_{2m}, y_1,\ldots ,y_k)$ such that $E$ is generated by $\partial_{y_1},\ldots ,\partial_{y_k}$. In particular, $(dx_1,\ldots ,dx_{2m},\alpha_1,\ldots ,\alpha_k)$ form a frame of $T^*M$ and $\omega$ can be written as: $\omegabar\wedge \alpha_1\wedge \ldots \wedge \alpha_k$ , where $\omegabar=\sum f_{ij}(x,y)dx_i\wedge dx_j$.
	We now can split the de Rham differential into three components $d=R + \diff^x + \diff^\alpha$, of bidegrees $(2,-1), (1,0), (0,1)$ with respect to the splitting $TM=D(\omega)\oplus E$, as in e.g. \cite{lychaginrubtsov94} (There is no $(-1,2)$ component, since $E$ is involutive).\\

	With respect to this decomposition $\omega$ is a $(2,k)$-form, the wedge product of the $(0,k)$-form $\alpha_1\wedge \ldots \wedge \alpha_k$ with the  $(2,0)$-form $\omegabar$. In particular,  $d\omega=0$ implies $R\omega=0$, i.e.
	$$
		0= \sum_{i=1}^k\pm R(\alpha_i)\wedge \omegabar \wedge \alpha_1\wedge \ldots \wedge \widehat\alpha_i\wedge \ldots \wedge \alpha_k.
	$$
	Since the different $\alpha_1\wedge \ldots \wedge \widehat\alpha_i\wedge \ldots \wedge \alpha_k$  (bidegree $(0,k-1)$) terms can not cancel, this means that $R(\alpha_i)\wedge \omegabar$ has to be zero for each $i$. The $R(\alpha_i)$ are $(2,0)$ -forms, i.e., $y$-dependent 2-forms in the $x$ coordinates and $\omegabar$ is a ($y$-dependent) symplectic form in the $x$-directions. \\

	Lepage's decomposition theorem for symplectic forms \cite[Section 16]{libermannmarle1987}, says that the operator $\wedge \omegabar$ from $\Omega^{(l,0)}(M)$ to $\Omega^{(l+2,0)}(M)$ is injective for $l\leq m-1$. Since we assume $m>2$, it is injective for  $l=2$ implying that $R(\alpha_i)=0$ for all $i$. This means that the ideal spanned by $\alpha_1,\ldots ,\alpha_k$ is differential, i.e., that $D(\omega)$ is involutive.
\end{proof}

\begin{cor}[Theorem 10.1 in \cite{MR801210}]
	A decomposable multisymplectic 3-form $\omega$ in a manifold of dimension $>5$ is always flat. Decomposable means that it can be written as the wedge product of a 1-form with a 2-form (i.e., that $F(\omega)$ is 1-dimensional). In \cite{MR801210}, this corollary is proven under the additional assumption of infinitesimal transitivity of the form, i.e., the case when any tangent vector can be locally extended to a vector field preserving $\omega$.
\end{cor}

\section{Symmetries}
\label{sec:symmetry}
In this section, we discuss the symmetry-behavior of dvs-forms. We restrict our attention to the case of involutive $D(\omega)$, i.e., the flat case. An investigation for non-flat case for 3-forms in 5-dimensional space can be found in \cite{MR801210}.\\

\subsection{Global symmetries}
Before turning to infinitesimal symmetries, let us briefly note that there is an immediate difference between symplectic forms and higher degree dvs-forms:

\begin{lem}
	Let $M$ be a connected manifold and $\omega$ a dvs-form. If $\mathrm{rank}(F(\omega))=0$, then the group $\mathrm{Diff}(M,\omega)$ acts $n$-transitively on $M$ for all $n\in\mathbb N$.\footnote{The same holds for volume forms, i.e., when $\mathrm{rank}(F(\omega))=\mathrm{dim}(M)$. However, strictly speaking, this is not a density-valued symplectic form in our sense.} When $0<\mathrm{rank}(F(\omega))<\mathrm{dim}(M)$, the group does not act 2-transitively.
\end{lem}

Here, we say that an action is called $n$-transitive, if for any two tuples of distinct points $(p_1,\ldots ,p_n)$, $(q_1,\ldots ,q_n)$ in $M$, there is an element $\phi\in \mathrm{Diff}(M,\omega)$ such that $\psi(p_i)=q_i~\forall i$. In the symplectic case (and for volume forms), the $n$-transitivity has been proven in \cite{MR0236961}. In all the other cases, the involutive distribution $D(\omega)$ induces a foliation preserved by $\mathrm{Diff}(M,\omega)$. Hence, no element of $\mathrm{Diff}(M,\omega)$ can map a pair of points on the same leaf of the foliation to two points on different leaves of the foliation.\\

For the standard form \eqref{vardarbouxform} on $M=\mathbb R^{k+2m}$ the group $\mathrm{Diff}(M,\omega)$ acts 1-transitively. However, for $k,m\neq 0$ we can easily construct examples where $\mathrm{Diff}(M,\omega)$ does not even act 1-transitively by removing 1 point from $M$, hence rendering one leaf non-diffeomorphic to the others.\\

However, for all dvs-forms, the group $\mathrm{Diff}(M,\omega)$ stays very big. Even if we consider only diffeomorphisms fixing the leaves of the foliation, there is still an infinite-dimensional group of symplectic diffeomorphisms acting on each leaf. We make this more precise in the next section.

\subsection{Infinitesimal Symmetries}
In this subsection, we are interested in the vector fields preserving $\omega$, i.e., the space $\mathfrak X(M,\omega)=\{X~|~L_X\omega=0\}$.\\

Let us look at the model form $\omega=\omegabar\wedge \nu\in \Omega^{k+2}(\mathbb R^{2m+k})$, where $\nu:=dy_1\wedge \ldots  \wedge dy_k$ and $\omegabar=\sum_{i=1}^m \diff x_{2i-1}\wedge \diff x_{2i}$. Let $V\in \mathfrak X(M)$ be a vector field. We denote by $V(x,y)=V^x(x,y)+V^y(x,y)$ its splitting onto its  $x$-directions and $y$-directions respectively. We first calculate
$$
	L_V\omega=   L_V\omegabar\wedge \nu+\omegabar\wedge(L_V\nu)  = d\iota_{V_x}\omegabar\wedge\nu + \omegabar\wedge (d\iota_{V^y}\nu).
$$
Now, we again use the fact that, in our case, the de Rham complex is bigraded and split $d$ into $d^x+d^y$.

$$
	L_V\omega= \omegabar\wedge (d^x\iota_{V^y}\nu) + d^x\iota_{V_x}\omegabar\wedge\nu +\omegabar \wedge (d^y\iota_{V^y}\nu)+ d^y\iota_{V_x}\omegabar\wedge\nu.
$$

Since there are no non-trivial $(1,k+1)$-forms the $d^y\iota_{V_x}\omegabar \wedge \nu$-component vanishes. Furthermore, $d^y\iota_{V^y}\nu=L^y_{V^y}\nu=div^y(V^y)\nu$, where the $y$-superscript in the Lie derivative and divergence indicates that we treat the $x$-variables as constants. Similarly $ d^x\iota_{V_x}\omegabar$ can be rewritten to $ L^x_{V_x}\omegabar$  The above equation hence reads:

$$L_V\omega= \omegabar\wedge (d^x\iota_{V^y}\nu) + L^x_{V_x}\omegabar\wedge \nu +\omegabar\wedge div^y({V^y})\nu.$$

Due to the bigrading, the equation $L_V\omega=0$ translates to the equations:
\begin{align*}
	 & 0=\omegabar\wedge(d^x\iota_{V^y}\nu)                    \\
	 & 0=(L^x_{V_x}\omegabar+div^y({V^y}) \omegabar)\wedge \nu
\end{align*}

Again due to Lepage's theorem, which we already used in Theorem \ref{thm:autoflat}, wedging $(1,k-1)$-forms with $\omegabar$ is an injective operation, such that the first equation simplifies to $d^x\iota_{V^y}\nu=0$, which means that $V^y$ can not depend on $x$, i.e., $V^y(x,y)=V^y(y)$. Similarly, the second equation reduces to $$L^x_{V_x}\omegabar=-div^y({V^y}) \omegabar.$$

Since we have already established that $V^y$ does not depend on $x$, we can consider this equation for each value of $y$ separately. There, this equation becomes the equation for $V^x$ to be a conformally symplectic vector field with conformal constant $-div^y({V^y})$. Conformally symplectic vector fields have been studied in \cite{confhamsys}. When the conformal constant is zero, those are just symplectic vector fields. Otherwise, we have the following:
\begin{prop}[\cite{confhamsys}]\label{prop:conf}
	Let $(P,\eta)$ be a symplectic manifold and $c\neq 0$ a constant.
	\begin{enumerate}
		\item Conformal Hamiltonian vector fields for $c$ (i.e., vector fields $X$ satisfying $L_X\eta=c\eta$) exist if and only if $\eta$ is exact.
		\item Let $d\theta=\eta$ and $\mathcal E$ the unique vector field satisfying $\iota_{\mathcal E}\eta=\theta$. Then any conformally symplectic vector field for $c$ is of the form $X=Y+c\mathcal E$, where $Y$ is a symplectic vector field.
		\item Hence, if the de Rham cohomology group $H^1(P)$ is trivial, then any conformally symplectic vector field for $c$ can be written as $X_H+c\mathcal E$, where $X_H$ denotes the Hamiltonian vector field of some $H\in C^{\infty}(P)$.
	\end{enumerate}
\end{prop}
In our case, we can choose a $d^x$-primitive of $\omegabar$, e.g. $\theta=x_1\wedge dx_2+\ldots +x_{2m-1}\wedge dx_{2m}$ and get a corresponding vector field $\mathcal E$. We can furthermore construct for any smooth function $H(x,y)$ a Hamiltonian vector field $X^x_H$ by solving $d^xH=\iota_{X^x_H}\omegabar$. In total, we have shown:

\begin{thm} Let $M=\mathbb R^{2m}\times \mathbb R^k$, $\omega=\omegabar\wedge \nu$ and $\mathcal E$ be as above. Then
	$$\mathfrak X(M,\omega)=\{X^x_H- div^y{Y}\mathcal E + Y~|~Y\in \mathfrak X(\mathbb R^k),~H\in C^{\infty}(\mathbb R^{2m}\times \mathbb R^k)
		\}.$$
\end{thm}

\begin{rem}The above argumentation can be carried out with any product of an exact symplectic manifold $(P,\omegabar=d\theta)$ with a manifold with volume form $(N,\nu)$.
	In the absence of $H^1(P)=0$, we just have to substitute the parameter-dependent Hamiltonian vector field with a parameter-dependent symplectic one:
	$$\mathfrak X(P\times N,\omegabar\wedge\nu )=\{X- div^y{Y}\mathcal E + Y~|~Y\in \mathfrak X(N),~X:P\times N\to TP \mathrm{~s.th.~} X(\cdot,y)\in\mathfrak X(P,\omegabar)\forall y\in N
		\}.$$

	However, if $\omegabar$ is non-exact (e.g. $P$ compact), then by Proposition \ref{prop:conf} there are no conformally symplectic vector fields and we get
	$$\mathfrak X(P\times N,\omegabar\wedge\nu )=\{X + Y~|~Y\in \mathfrak X(N,\nu),~X:P\times N\to TP \mathrm{~s.th.~} X(\cdot,y)\in\mathfrak X(P,\omegabar)\forall y\in N
		\}.$$
	This means that not every germ of a Hamiltonian vector field can be extended to a global Hamiltonian vector field.
\end{rem}

\begin{rem}
	The notion of Hamiltonian vector fields extends to multisymplectic geometry. A vector field on a multisymplectic manifold $(M,\omega)$ is a vector field $X$ such that $\iota_X\omega=d\alpha$ for some differential form $\alpha$, which is then called the Hamiltonian form for $X$. In the case of $M=P\times N$ and $\omega=\omegabar\wedge \nu$ with $\omegabar=d\theta$ as above, there is a very natural candidate for a Hamiltonian form of $X^x_H- div^y{Y}\mathcal E + Y$, it is given by $\alpha=H\nu+\iota_Y\nu\wedge\theta$.
\end{rem}

\begin{exmp} In the case $M=\mathbb R^{2m+1}$ (i.e., $k=1$), we recover the description of \cite{MR801210}, that symmetries of $(M,\omega)$ are parametrized by two functions $R=R(y)$ (corresponding to the vector field $R\partial_y$) and $H=H(x,y)$.
\end{exmp}

\section{Relation to cosymplectic geometry}
\label{sec:cosymplectic}

Cosymplectic structures have been introduced in \cite{libermann}, cf. also \cite{cosympsurvey} for a recent survey. They can be seen as an odd counterpart to symplectic structures and naturally appear on the singular loci of $b$-symplectic manifolds \cite[Theorem 3.2]{logsympyieldscosymp}.

\begin{defn} Let $M$ be a $(2n+1)$-dimensional manifold with $n>1$. A cosymplectic structure on $M$ is given by a pair $(\alpha,\beta)\in \Omega^1(M)\times \Omega^2(M)$ such that:
	\begin{itemize}
		\item $d\alpha=0$ and $d\beta=0$.
		\item $\alpha\wedge \beta^n$ is non-degenerate (i.e., a volume form).
	\end{itemize}
\end{defn}

By construction, for any cosymplectic manifold, the form $\beta\wedge \alpha$ is a dvs-form. (This has been considered in the context of reduction in \cite{coredgroupoids}). Since $\alpha$ is closed, $D(\omega)=\ker(\alpha)$ is involutive, hence Theorem \ref{thm:ourdarboux} is applicable. Moreover, we can choose coordinates  $(x_1,\ldots ,x_{2n},y)$ adapted to both the foliation $D(\omega)$ and the complementary one-dimensional foliation $E=\{v| ~ \iota_v\beta=0\}$. These coordinates satisfy $ \iota_{\partial_{x_i}}\alpha = 0$ and $\iota_{\partial_{y}}\beta = 0$. Upon rescaling $y$, we can even obtain $\iota_{\partial_{y}}\alpha = 1$, i.e., $\alpha=dy$. Correspondingly, since $d^y\beta=0$,  $\beta$ is completely independent of $y$ and can be flattened with the classical symplectic Darboux theorem. We recover the cosymplectic Darboux theorem:

\begin{thm}[\cite{CODARBOUX}] Given a cosymplectic structure, there exist local coordinates such that $\alpha=dy$ and $\beta=\sum_{i=1}^ndx_{2i-1}\wedge dx_{2i}$.
\end{thm}

There are several types of (infinitesimal) symmetries one can consider in the cosymplectic setting \cite{cosympham}:
\begin{itemize}
	\item $X$ is called cosymplectic if $L_X\alpha=0$ and $L_X\beta=0$.
	\item $X$ is called weakly co-Hamiltonian if it is cosymplectic and if $\iota_X\beta+\alpha(X)\alpha$ is exact.
	\item $X$ is called co-Hamiltonian if it is weakly co-Hamiltonian and $\alpha(X)=0$.
\end{itemize}
Cosymplectic vector fields, of course, automatically preserve the form $\omega=\beta\wedge \alpha$, i.e., are multisymplectic. Co-Hamiltonian vector fields are Hamiltonian (in the multisymplectic sense): If $\alpha(X)=0$, then a primitive $\gamma$ of $\iota_X\beta+\alpha(X)\alpha=\iota_X\beta$ induces the primitive $\gamma\wedge \alpha$ of $\iota_X(\beta\wedge \alpha)$. On the other hand, weakly co-Hamiltonian fields need not be (multisymplectically) Hamiltonian: If we take the standard Fubini-Study symplectic form $\beta$ on $\mathbb CP^2$ and consider its product with $\mathbb R$ (with the form $dy$), we obtain a cosymplectic manifold. The vector field $\partial_y$ is weakly co-Hamiltonian $\iota_{\partial_y}\beta+dy(\partial_y)dy=dy$ is exact. However, it is not Hamiltonian in the multisymplectic sense, since $\iota_{\partial_y}(\beta\wedge dy)=\beta$ is not exact.

\begin{rem}
	We finish by noting that $k$-cosymplectic forms in the sense of \cite[Definition 3.2.20]{thesisdefininingkcosymp} also provide instances of dvs-forms \footnote{A competing notion of $k$-cosymplectic of a different flavor not yielding a dvs-form in general was presented in \cite[Definition 4.1]{closedarbouxsurvey}}. Moreover, elliptic symplectic structures with zero residue also induce dvs-forms, even when their 2-form is not closed, i.e., when they do not induce 2-cosymplectic structures \cite[Proposition 5.2.27 and Equation (5.2.8)]{Witte2021}.
\end{rem}

\bibliographystyle{plain}
\bibliography{ref}

\begin{thebibliography}{10}

\bibitem{CODARBOUX}
Claude Albert.
\newblock Le th{\'e}or{\`e}me de r{\'e}duction de {Marsden}-{Weinstein} en
  g{\'e}om{\'e}trie cosymplectique et de contact. ({The} {Marsden}-{Weinstein}
  reduction theorem in cosymplectic and contact geometry).
\newblock {\em J. Geom. Phys.}, 6(4):627--649, 1989.

\bibitem{MR0236961}
William~M. Boothby.
\newblock Transitivity of the automorphisms of certain geometric structures.
\newblock {\em Trans. Amer. Math. Soc.}, 137:93--100, 1969.

\bibitem{cosympsurvey}
Beniamino Cappelletti-Montano, Antonio De~Nicola, and Ivan Yudin.
\newblock A survey on cosymplectic geometry.
\newblock {\em Rev. Math. Phys.}, 25(10):55, 2013.
\newblock Id/No 1343002.

\bibitem{logsympyieldscosymp}
Gil~R. Cavalcanti.
\newblock Examples and counter-examples of log-symplectic manifolds.
\newblock {\em J. Topol.}, 10(1):1--21, 2017.

\bibitem{coredgroupoids}
Daniel~López Garcia and Nicolas~Martinez Alba.
\newblock Reduction of cosymplectic groupoids, 2024.
\newblock arXiv:2403.03178.

\bibitem{closedarbouxsurvey}
Xavier Gràcia, Javier de~Lucas, Xavier Rivas, and Narciso Román-Roy.
\newblock On darboux theorems for geometric structures induced by closed forms,
  2023.
\newblock arXiv:2306.08556.

\bibitem{MR2882772}
Fr\'ed\'eric H{\'e}lein.
\newblock Multisymplectic formalism and the covariant phase space.
\newblock In {\em Variational problems in differential geometry}, volume 394 of
  {\em London Math. Soc. Lecture Note Ser.}, pages 94--126. Cambridge Univ.
  Press, Cambridge, 2012.

\bibitem{MR0334772}
Jerzy Kijowski.
\newblock A finite-dimensional canonical formalism in the classical field
  theory.
\newblock {\em Comm. Math. Phys.}, 30:99--128, 1973.

\bibitem{krieglMichor}
Andreas Kriegl and Peter~W. Michor.
\newblock {\em The convenient setting of global analysis}, volume~53 of {\em
  Math. Surv. Monogr.}
\newblock Providence, RI: American Mathematical Society, 1997.

\bibitem{Leski2022}
Laura Leski.
\newblock Flatness and non-flatness of multisymplectic manifolds, 2022.
\newblock Bachelor's Thesis.

\bibitem{libermann}
Paulette Libermann.
\newblock Sur les automorphismes infinit{\'e}simaux des structures
  symplectiques et des structures de contact.
\newblock Centre {Belge} {Rech}. {Math}., {Colloque} {G{\'e}om}. {Diff{\'e}r}.
  {Globale}, {Bruxelles} du 19 au 22 {D{\'e}c}. 1958, 37-59 (1959)., 1959.

\bibitem{libermannmarle1987}
Paulette Libermann and Charles-Michel Marle.
\newblock {\em Symplectic geometry and analytical mechanics. {Transl}. from the
  {French} by {Bertram} {Eugene} {Schwarzbach}}, volume~35 of {\em Math. Appl.,
  Dordr.}
\newblock Springer, Dordrecht, 1987.

\bibitem{lychaginrubtsov94}
Valentin Lychagin and Vladimir Rubtsov.
\newblock Non-holonomic filtration: {Algebraic} and geometric aspects of
  non-integrability.
\newblock In {\em Geometry in partial differential equations}, pages 189--214.
  Singapore: World Scientific, 1994.

\bibitem{martin}
Geoffrey Martin.
\newblock A {Darboux} theorem for multi-symplectic manifolds.
\newblock {\em Lett. Math. Phys.}, 16(2):133--138, 1988.

\bibitem{martinet70}
Jean Martinet.
\newblock Sur les singularit\'es des formes diff\'erentielles.
\newblock {\em Ann. Inst. Fourier (Grenoble)}, 20(fasc. 1):95--178, 1970.

\bibitem{confhamsys}
Robert McLachlan and Matthew Perlmutter.
\newblock Conformal {Hamiltonian} systems.
\newblock {\em J. Geom. Phys.}, 39(4):276--300, 2001.

\bibitem{thesisdefininingkcosymp}
B.~{Osorno Torres}.
\newblock {\em Codimension-one Symplectic Foliations: Constructions and
  Examples}.
\newblock Doctoral thesis, Universiteit Utrecht, September 2015.

\bibitem{MR2892558}
Christopher~L. Rogers.
\newblock {$L_\infty$}-algebras from multisymplectic geometry.
\newblock {\em Lett. Math. Phys.}, 100(1):29--50, 2012.

\bibitem{RYVKIN20199}
Leonid Ryvkin and Tilmann Wurzbacher.
\newblock An invitation to multisymplectic geometry.
\newblock {\em Journal of Geometry and Physics}, 142:9--36, 2019.

\bibitem{sevwurz}
Gabriel Sevestre and Tilmann Wurzbacher.
\newblock Lagrangian submanifolds of standard multisymplectic manifolds.
\newblock In {\em Geometric and harmonic analysis on homogeneous spaces.
  Selected papers of the 5th Tunisian-Japanese conference, PTJC 2017, Mahdia,
  Tunisia, December 17--21, 2017}, pages 191--205. Cham: Springer, 2019.

\bibitem{cosympham}
Stephane Tchuiaga, Franck Houenou, and Pierre Bikorimana.
\newblock On cosymplectic dynamics. {I}.
\newblock {\em Complex Manifolds}, 9:114--137, 2022.

\bibitem{MR801210}
J.~F. Turiel.
\newblock {\em Classification locale des 3-formes ferm\'ees
  infinit\'esimalement transitives \`a cinq variables}, volume~30 of {\em
  Cahiers Math\'ematiques Montpellier [Montpellier Mathematical Reports]}.
\newblock Universit\'e des Sciences et Techniques du Languedoc, U.E.R. de
  Math\'ematiques, Montpellier, 1984.

\bibitem{vaisman2016cohomology}
I.~Vaisman.
\newblock {\em Cohomology and Differential Forms}.
\newblock Dover Books on Mathematics. Dover Publications, 2016.

\bibitem{vanzuraCharacterizationOneType2004}
Ji{\v r}{\'i} Van{\v z}ura.
\newblock Characterization of one type of multisymplectic 3-forms in odd
  dimensions.
\newblock In {\em The Proceedings of the 23rd Winter School ``{{Geometry}} and
  Physics'', {{Srn{\'i}}}, {{Czech Republic}}, {{January}} 18--25, 2003}, pages
  203--209. Palermo: Circolo Matem{\`a}tico di Palermo, 2004.

\bibitem{ww}
Maxime Wagner and Tilmann Wurzbacher.
\newblock Equivalent formulations of hamiltonian dynamics on multi cotangent
  bundles, 2024.
\newblock arXiv:2410.21068.

\bibitem{Witte2021}
Aldo Witte.
\newblock {\em Between generalized complex and Poisson geometry}.
\newblock PhD thesis, Utrecht University, 2021.

\end{thebibliography}

\end{document}